\documentclass[12pt]{article}
\usepackage{jhep-mod}
\usepackage{bm}
\usepackage{amssymb,amsmath,amsthm}
\usepackage{mathrsfs}
\usepackage[utf8]{inputenc}
\usepackage{enumerate}
\hypersetup{colorlinks=true} 
\usepackage{appendix}
\usepackage{graphicx}
\usepackage{dcolumn}
\usepackage{bm}
\usepackage{multirow}
\usepackage{float}
\usepackage{tikz}
\usepackage[normalem]{ulem}
\definecolor{purple}{rgb}{1,0,1}
\definecolor{lime}{HTML}{A6CE39} % needs xcolor

\newcommand{\blue}[1]{{\color{blue} #1}}

\newtheorem{sufficient}{Sufficient condition}
\definecolor{lime}{HTML}{A6CE39}
\newcommand{\orcidicon}{%
	\begin{tikzpicture}
	\draw[lime, fill=lime] (0,0) 
		circle [radius=0.16] 
		node[white] {{\fontfamily{qag}\selectfont \tiny ID}};
	\draw[white, fill=white] (-0.0625,0.095) 
		circle [radius=0.007];
	\end{tikzpicture}
	\hspace{-2mm}
}
\newcommand\orcidMatt{{\href{https://orcid.org/0000-0003-1088-6485}{\orcidicon}}}
\begin{document}
%========================================================

\title{\huge {Verifying the Firoozbakht, Nicholson, and Farhadian conjectures up to the 81$^{st}$ maximal prime gap.}}

\author{\Large Matt Visser\orcidMatt{}}
%========================================================
%========================================================
%========================================================
%========================================================
\affiliation{School of Mathematics and Statistics, Victoria University of Wellington, \\
PO Box 600, Wellington 6140, New Zealand}
%========================================================
%========================================================
\emailAdd{matt.visser@sms.vuw.ac.nz}
%========================================================
%========================================================

\abstract{
\parindent0pt
\parskip7pt
The Firoozbakht, Nicholoson, and Farhadian conjectures can be phrased in terms of increasingly powerful conjectured bounds on the prime gaps $g_n := p_{n+1}-p_n$. 
\[
g_n  \leq p_n \left(p_n^{1/n} -1 \right)\qquad\qquad\qquad (n \geq 1; \; Firoozbakht).
\]
\[
g_n  \leq p_n \left((n\ln n)^{1/n} -1 \right)\qquad\qquad (n>4; \; Nicholson).
\]
\[
g_n  \leq p_n \left( \left(p_n {\ln n\over\ln p_n}\right)^{1/n} -1 \right)\qquad (n>4; \; Farhadian).
\]

While a general proof of any of these conjectures is far out of reach I shall show that all three of these conjectures are unconditionally and explicitly verified for all primes below the location of the 81$^{st}$ maximal prime gap, certainly for all primes $p <2^{64}$. 
For the Firoozbakht conjecture this is a very minor improvement on currently known results, for the Nicholson and Farhadian conjectures this may be more interesting.
 
\bigskip
{\sc Date:} 1 April 2019; 8 April 2019; \LaTeX-ed \today

\bigskip
{\sc ArXiv:} 1904.00499

\medskip
{\sc Keywords:} 
primes; prime gaps; Firoozbakht conjecture; Nicholson conjecture; \\
Farhadian conjecture.

\medskip
{\sc MSC:} 
11A41 (Primes);  11N05 (Distribution of primes).

\medskip
{\sc Sequences:} 
A005250 A002386 A005669 A000101 A107578  A246777 A246776
\vfill
}

\notoc
\maketitle

%========================================================
\def\tr{{\mathrm{tr}}}
\def\diag{{\mathrm{diag}}}
\parindent0pt
\parskip7pt

\null
\vspace{-50pt}
%========================================================
\section{Introduction}
%========================================================
%---------------------------------------------------------------------------------------------------------------------------------------------
\label{S:intro}
%---------------------------------------------------------------------------------------------------------------------------------------------
\def\ceil{{\mathrm{ceil}}}
\def\floor{{\mathrm{floor}}}
\def\frac{{\mathrm{frac}}}
\def\N{{\mathbb{N}}}
\def\Z{{\mathbb{Z}}}
\def\P{{\mathbb{P}}}
\def\implies{\Longrightarrow}
\newtheorem{conjecture}{Conjecture}{}
\parindent0pt
\parskip7pt

The Firoozbakht, Nicholson, and Farhadian conjectures would, if proved to be true, impose increasingly strong constraints on the distribution of the primes;  this distribution being a fascinating topic that continues to provide many subtle and significant open questions~\cite{Ribenboim:91, Ribenboim:96, Ribenboim:04, Wells, Cramer:1919, Cramer:1936, Goldston, Rosser:38, Rosser:41, Cesaro, Cippola, Rosser:62, Sandor, Dusart-1999, Lowry-Duda, Dusart:2010, Trudigan:2014, Dusart2, Axler:2017, Visser:Lambert, Andrica, Visser:Andrica, Visser:strong-Andrica}.
The Firoozbakht conjecture~\cite{Firoozbakht, Firoozbakht2, Kourbatov1, Kourbatov2, Kourbatov3} is normally phrased as follows.
\begin{conjecture}{(Firoozbakht conjecture, two most common versions)}
\begin{equation}
(p_{n+1})^{1\over n+1} \leq (p_{n})^{1\over n}  ; \qquad   \hbox{equivalently} \qquad    {\ln p_{n+1}\over n+1} \leq {\ln p_n\over n}; \qquad (n\geq 1).
\end{equation}
\end{conjecture}
To see why this conjecture might be somewhat plausible, use the standard inequalities $n \ln n < p_n < n \ln p_n$, which hold for $n\geq1$ and $n\geq 4$ respectively, and observe that
\begin{equation}
{\ln (n\ln n)\over n}  \leq {\ln p_n\over n} \leq {\ln^2 p_n\over p_n};  \qquad (n\geq 1; n \geq 4).
\end{equation}
Now ${\ln (n\ln n)\over n}$ is monotone decreasing for $n\geq 5$, and ${\ln^2 p_n\over p_n}$ is monotone decreasing for $p_n> 7$. So  for $n\geq 5$, corresponding to $p_n\geq 11$, the function ${\ln p_n\over n}$ is certainly bounded between two monotone decreasing functions; the overall trend is monotone decreasing. The stronger conjecture that  ${\ln p_n\over n}$ is itself monotone decreasing depends on fluctuations in the distribution of the primes $p_n$; fluctuations which can be rephrased in terms of the prime gaps $g_n := p_{n+1}-p_n$. 

\enlargethispage{30pt}
Indeed, Kourbatov~\cite{Kourbatov1} using results on first occurrence prime gaps
has recently verified Firoozbakht's conjecture to hold for all primes $p<4\times 10^{18}$. 
Furthermore Kourbatov~\cite{Kourbatov2} has also derived a \emph{sufficient} condition for the Firoozbakht conjecture to hold:
\begin{equation}
g_n \leq \ln^2 p_n - \ln p_n - 1.17;  \qquad (n\geq10; \; p_n\geq29).
\end{equation}
Using tables of  first occurrence prime gaps and maximal prime gaps Kourbatov has now extended this discussion~\cite{Kourbatov3}, and subsequently verified that Firoozbakht's conjecture holds for all primes $p < 1\times 10^{19}$. More recently (2018), two additional maximal prime gaps have been found~\cite{g80}, so that Kourbatov's  arguments now certainly verify the Firoozbakht conjecture up to the 80$^{th}$ maximal prime gap --- more precisely, for all primes below currently unknown location of the 81$^{st}$ maximal prime gap ---  {though we do now know (September 2018) that $p^*_{81}> 2^{64}$~\cite{Nicely}, see also~\cite{Kourbatov3}.  So certainly the Firoozbakht conjecture holds for all primes $p < 2^{64} = 18,446,744,073,709,551,616 \approx 1.844\times 10^{19}$. Note that this automatically verifies a strong form of Cram\'er's conjecture
\begin{equation}
g_n \leq \ln^2 p_n;  \qquad (n\geq5; \; p_n\geq11),
\end{equation}
at least for all primes $p <2^{64}\approx1.844\times 10^{19}$.

What is trickier with Kourbatov's techniques is to say anything useful about the slightly stronger Nicholson~\cite{Nicholson} and Farhadian~\cite{Farhadian1,Farhadian2} conjectures, and it is this issue we shall address below.

%========================================================
\section{Firoozbakht, Nicholson, and Farhadian}
%========================================================
When comparing the Firoozbakht conjecture with the slightly stronger Nicholson and Farhadian conjectures it is useful to work with the ratio of successive primes, $p_{n+1}/p_n$.
\begin{conjecture}{(Firoozbakht/Nicholson/Farhadian conjectures; successive primes)}
\begin{align}
{(p_{n+1}/p_n)^n}& \leq p_n\quad &(n& \geq 1; \; Firoozbakht).
\\
{(p_{n+1}/p_n)^n}& \leq n \ln n\quad &(n&>4; \; Nicholson).
\\
{(p_{n+1}/p_n)^n}& \leq p_n {\ln n\over\ln p_n} \quad &(n&>4; \; Farhadian).
\end{align}
\end{conjecture}
When phrased in this way the standard inequalities $ n \ln n < p_n < n \ln p_n$ show that Farhadian $\implies$ Nicholson $\implies$ Firoozbakht. 
To study the numerical evidence in favour of  these conjectures it is useful to convert them into statements about the prime gaps $g_n := p_{n+1}-p_n$. 
\begin{conjecture}{(Firoozbakht/Nicholson/Farhadian conjectures; prime gap version)}
\begin{align}
&g_n  \leq p_n \left(p_n^{1/n} -1 \right)\quad &(n& \geq 1; \; Firoozbakht).
\\
&g_n  \leq p_n \left((n\ln n)^{1/n} -1 \right)\quad &(n&>4; \; Nicholson).
\\
&g_n  \leq p_n \left( \left(p_n {\ln n\over\ln p_n}\right)^{1/n} -1 \right)\quad &(n&>4; \; Farhadian).
\end{align}
This can further be rephrased as:
\begin{align}
&g_n  \leq p_n \left(\exp\left({\ln p_n\over n}\right) -1 \right)\quad &(n& \geq 1; \; Firoozbakht).
\\
&g_n  \leq p_n \left(\exp\left({\ln (n\ln n)\over n}\right)  -1 \right)\quad &(n&>4; \; Nicholson).
\\
&g_n  \leq p_n \left( \exp\left({1\over n} \ln \left(p_n {\ln n\over\ln p_n}\right)\right) -1 \right)\quad &(&n>4; \; Farhadian).
\end{align}
\end{conjecture}
These inequalities are all of the form $g_n\leq f(p_n, n)$, with $f(p_n,n)$ a function of both $p_n$ and $n$. 

While  $p_n$ and $n$ are both monotone increasing, unfortunately $f(p_n,n)$ is not guaranteed to be monotone increasing, so one would have to check each individual $n$ independently. 
So our strategy will be to seek to find suitable \emph{sufficient} conditions for the Firoozbakht/Nicholson/Farhadian conjectures of the form $g_n \leq f(n)$, with the function $f(n)$ being some monotone function of its argument. Once this has been achieved we can develop an argument using maximal prime gaps. 

%========================================================
\section{\leftline{Sufficient condition for the Nicholson/Firoozbakht conjectures}}
%========================================================

Using the fact that $e^x-1> x$ we deduce a sufficient condition for the Nicholson conjecture (which is then automatically also sufficient for the Firoozbakht conjecture).
\begin{sufficient}{(Nicholson/Firoozbakht)}
\begin{equation}
g_n < {p_n \ln (n\ln n)\over n}; \qquad\qquad (n>4; \;\; n \geq 1).
\end{equation}
\end{sufficient}
Now use Dusart's result~\cite{Dusart-1999} that for $n\geq 2$ we have $p_n > n(\ln(n\ln n) - 1)$ to deduce the stronger sufficient condition
\begin{sufficient}{(Nicholson/Firoozbakht)}
\begin{equation}
g_n < f(n) = (\ln(n\ln n) - 1) \ln (n\ln n); \qquad (n>4; \; n \geq 2).
\end{equation}
\end{sufficient}
\emph{A posteriori} we shall verify that this last condition is strong enough to be useful, and weak enough to be true over the domain of interest.

%========================================================
\section{{Verifying the Firoozbakht and Nicholson conjectures \\ for all primes {$p < 2^{64}$}}}
%========================================================
%---------------------------------------------------------------------------------------------------------------------------------------------
\label{S:verify-firoozbakht}
%---------------------------------------------------------------------------------------------------------------------------------------------
This is a variant of the argument given for the Andrica conjecture in references~\cite{Visser:Andrica,Visser:strong-Andrica}.
Consider the maximal prime gaps: Following a minor modification of the notation of references~\cite{Visser:Andrica,Visser:strong-Andrica}, let the quartet $(i, g^*_i, p^*_i, n_i^*)$ denote the $i^{th}$ maximal prime gap; of width $g^*_i$, starting at the $n_i^*$th prime $p^*_i = p_{n_i^*}$. 
(See see the sequences A005250, A002386, A005669, A000101, A107578.)
As of April 2019, some 80 such maximal prime gaps are known~\cite{g80,g75,gaps}, up to $g^*_{80}=1550$ and 
\begin{equation}
p^*_{80}= 18,361,375,334,787,046,697 >1.836\times 10^{19},
\end{equation}
which occurs at 
\begin{equation}
n^*_{80} = 423,731,791,997,205,041 \approx 423\times10^{15}. 
\end{equation}
One now considers the interval $[p^*_i, p^*_{i+1}-1]$, from the lower end of the $i^{th}$ maximal prime gap to just below the beginning of the $(i+1)^{th}$ maximal prime gap. Then everywhere in this interval
\begin{equation}
\forall p_n\in [p^*_i,p^*_{i+1}-1]  \qquad   g_n \leq g^*_i; \qquad  f(n_i^*) \leq  f(n).
\end{equation}
Therefore, if the sufficient condition for the Nicholson/Firoozbakht conjectures holds at the beginning of the  interval $p_n\in [p^*_i,p^*_{i+1}-1]$, then it certainly holds on the entire interval. (Note that for the Nicholson/Firoozbakht conjectures, in addition to knowing the $p_i^*$, it is also essential to know all the $n_i^* = \pi(p_i^*)$ in order for this particular verification procedure to work; for the Andrica conjecture one can quietly discard the $n_i^* = \pi(p_i^*)$ and only work with the $p_i^*$~\cite{Visser:Andrica,Visser:strong-Andrica}.) 

%\clearpage
Explicitly checking a table of maximal prime gaps~\cite{g80,g75,gaps},  both of the Nicholson and Firoozbakht conjectures certainly hold on the interval $[p_5^*,p^*_{81}-1]$, from $p^*_5=89$  up to just before the beginning of the 81$^{st}$ maximal prime gap, $p^*_{81}-1$, even if we do not yet know the value of $p^*_{81}$. Then explicitly checking the primes below 89 the Firoozbakht conjecture holds for all primes $p$ less than $p_{81}^*$,
while the Nicholson conjecture holds for all primes $p$ less than $p_{81}^*$, except $p\in\{2,3,5,7\}$.
Since we do not explicitly know $p^*_{81}$, (though an exhaustive search has now verified that $p^*_{81}> 2^{64}$~\cite{Nicely}, see also~\cite{Kourbatov3}),  a safe fully explicit statement is that both the Firoozbakht and Nicholson conjectures are verified for all primes $p < 2^{64} \approx 1.844\times 10^{19}$.

%========================================================
\section{{Sufficient condition for the Farhadian conjecture}}
%========================================================

The Farhadian conjecture is a little trickier to deal with. Again using the fact that $e^x-1> x$ we can deduce a sufficient condition.
\begin{sufficient}{(Farhadian)}
\begin{equation}
g_n < {p_n \ln \left(p_n {\ln n\over\ln p_n}\right) \over n} =  {p_n \left(\ln p_n + \ln\ln n - \ln\ln p_n\right) \over n}; \qquad (n>4).
\end{equation}
\end{sufficient}
Now inside the brackets use the lower bound $p_n \geq n\ln n$ (valid for $n\geq1$), and the upper bound $p_n \leq n\ln (n\ln n) $ (valid for $n\geq6$).
This gives a new slightly stronger sufficient condition.
\begin{sufficient}{(Farhadian)}
\begin{equation}
g_n <   {p_n \left(\ln(n\ln n) + \ln\ln n - \ln\ln (n\ln(n\ln n) \right) \over n}; \qquad (n>6).
\end{equation}
\end{sufficient}
Now use Dusart's result~\cite{Dusart-1999} that for $n\geq 2$ we have $p_n > n(\ln(n\ln n) - 1)$ to deduce another yet even slightly stronger sufficient condition. 
\begin{sufficient}{(Farhadian)}
\begin{equation}
g_n < f(n) = (\ln(n\ln n) - 1) \left(\ln(n\ln n) + \ln\ln n - \ln\ln (n\ln(n\ln n) \right); \quad (n>6).
\end{equation}
\end{sufficient}
It is now a somewhat tedious exercise in elementary calculus to verify that this function $f(n)$ is indeed monotone increasing as a function of $n$. 
\emph{A posteriori} we shall verify that this last sufficient condition is strong enough to be useful, and weak enough to be true over the domain of interest.

%========================================================
\section{{Verifying the Farhadian conjecture for all primes $p<2^{64}$}}
%========================================================
%---------------------------------------------------------------------------------------------------------------------------------------------
\label{S:verify-farhadian}
%---------------------------------------------------------------------------------------------------------------------------------------------

The logic is the same as for the Firoozbakht and Nicholson conjectures. 
If the sufficient condition for the Farhadian  conjecture holds at the beginning of the  interval $p_n\in [p^*_i,p^*_{i+1}-1]$, then it certainly holds on the entire interval. Explicitly checking a table of maximal prime gaps~\cite{g80,g75,gaps},  the Farhadian conjecture certainly holds on the interval $[p_5^*,p^*_{81}-1]$, from $p^*_5=89$  up to just before the beginning of the 81$^{st}$ maximal prime gap, $p^*_{81}-1$, even if we do not yet know the value of $p^*_{81}$. Then explicitly checking the primes below $p_5^* = 89$ the Farhadian conjecture is verified to hold for all primes $p$ less than $p_{81}$, except $p\in\{2,3,5,7\}$.
Since we do not explicitly know $p^*_{81}$, 
(though an exhaustive search has now verified that $p^*_{81}> 2^{64}$~\cite{Nicely}, see also~\cite{Kourbatov3}),  a safe fully explicit statement is that the Farhadian conjecture is verified for all primes $p < 2^{64} \approx 1.844\times 10^{19}$.

%------------------------------------------------
\section{Discussion}\label{S:Discussion}
%------------------------------------------------
\def\li{{\mathrm{li}}}
While Kourbatov's recent work~\cite{Kourbatov1,Kourbatov2,Kourbatov3}  yields a useful and explicit domain of validity for the Firooz\-bakht conjecture, 
(ultimately, see~\cite{Nicely} and~\cite{Kourbatov3},  for all primes $p < 2^{64}$), the present article first slightly extends this domain of validity (all primes $p < p^*_{81}$), 
and second and more significantly obtains identical domains of validity for the related but slightly stronger Nicholson and Farhadian conjectures. 
The analysis has been presented in such a way that it can now be semi-automated. 

Upon discovery, every new maximal prime gap $g^*_i$ can, 
as long as one can also calculate the corresponding 
$n^*_i = \pi(p^*_i)$,  see for instance~\cite{LMO}, be used to push the domain
of validity a little further.    

Some cautionary comments are in order: Verification of these conjectures up to some maximal prime, 
however large, does not guarantee validity for all primes. Note that by the prime number theorem $\pi(n) \sim \li(n)$ so
\begin{equation}
{\ln(p_n)\over n} = {\ln(p_n)\over \pi(p_n)}  \sim {\ln(p_n)\over \li(p_n)}.
\end{equation}
Now certainly $\ln(p)/\li(p)$ is monotone decreasing, which is good. 
On the other hand $\pi(x)-\li(x)$ changes sign infinitely often, (this is the Skewes phenomenon~\cite{Skewes1,Skewes2,Skewes3,Skewes4}), so that the monotone decreasing function ${\ln(p_n)/ \li(p_n)}$ both over-estimates and under-estimates the quantity of interest ${\ln(p_n)/n}$, which is not so good. 
Now this observation does not disprove the Firoozbakht conjecture, but it does indicate where there might be some potential difficulty. 

On a more positive note, the Firoozbakht conjecture most certainly must hold when averaged over suitably long intervals. It is an elementary consequence of the Chebyshev theorems that $p_m p_n > p_{m+n}$, see \cite{Ribenboim:91, Ribenboim:96, Ribenboim:04}. But then $p_n^2 > p_{2n}$, and $p_n^3 > p_n p_{2n}> p_{3n}$. In general $(p_n)^m > p_{nm}$ and so $\ln p_n > \ln p_{nm}/m$. Consequently
\begin{equation}
{\ln(p_n)\over n} > {\ln(p_{nm})\over nm}. 
\end{equation}
This is much weaker than the usual Firoozbakht conjecture, but enjoys the merit of being unassailably true.

%------------------------------------------------
\acknowledgments{
This research was supported by the Marsden Fund, administered by the Royal Society of New Zealand.  
I particularly wish to thank Alexei Kourbatov for useful comments. 
}
%-------------------------------------------------

\clearpage
%========================================================

%========================================================
\end{document}